\documentclass[12pt,a4paper, oneside]{amsart}

\usepackage{amsmath, nicefrac, amsthm, verbatim, amsfonts, mathtools, amssymb, upgreek, xcolor, bbm}
\usepackage{graphics, xspace, enumerate}
\usepackage{stix}
\usepackage{dsfont}
\usepackage[a4paper,margin=2.5cm]{geometry}
\usepackage[bb=boondox]{mathalfa}

\usepackage{graphicx}
\usepackage[colorlinks=true,citecolor=red,urlcolor=blue,linkcolor=red,bookmarksopen=true,unicode=true,pdffitwindow=true]{hyperref}
\usepackage[english]{babel}
\usepackage[languagenames,fixlanguage]{babelbib}
\hypersetup{pdfauthor={}}
\hypersetup{pdftitle={}}

\usepackage{seqsplit,cleveref}

\theoremstyle{plain}
\newtheorem{theorem}{Theorem}[section]

\newtheorem*{theorem*}{Theorem}

\theoremstyle{definition}

\newtheorem{example}[theorem]{Example}

\newcommand {\Prob} {\ensuremath{\mathbb{P}}}
\newcommand {\R} {\ensuremath{\mathbb{R}}}
\newcommand {\ZZ} {\ensuremath{\mathbb{Z}}}

\newcommand{\df}{\coloneqq}

\newcommand{\E}{\mathrm{e}}

\newcommand{\Sp}{\mathsf{S}}

\newcommand{\T}{\mathbb{T}}

\newcommand{\D}{\mathrm{d}}

\title[Hoeffding's inequality for non-irreducible Markov models]{Hoeffding's inequality for non-irreducible Markov models}

\author[N.\ Sandri\'{c}]{Nikola Sandri\'{c}}
\address[Nikola\ Sandri\'{c}]{Department of Mathematics\\University of Zagreb\\ Zagreb\\Croatia}
\email{nikola.sandric@math.hr}

\author[S.\ \v Sebek]{Stjepan\ \v Sebek}
\address[Stjepan\  \v Sebek]{
	Department of Applied Mathematics\\
	Faculty of Electrical Engineering and Computing\\
	University of Zagreb\\ 
	Zagreb\\ 
	Croatia}
\email{stjepan.sebek@fer.hr}


\begin{document}
\allowdisplaybreaks[4]

\begin{abstract}
	In this article,  we establish  Hoeffding's inequality for  bounded  Lipschitz functions of  a class of not necessarily irreducible Markov models. The result complements the existing literature on this topic where Hoeffding's inequality for  bounded  measurable functions of  a class of  irreducible Markov models has been considered. Our  approach is based on the assumption of uniform ergodicity of the underlying Markov model in $\mathrm{L}^1$-Wasserstein space.
\end{abstract}

\maketitle


\section{Introduction}\label{S1}

Let $\{\xi_i\}_{i\in\ZZ_+}$ be a sequence of independent and bounded random variables (defined on a probability space $(\Omega,\mathcal{F}, \Prob)$), and let $S_t\df\xi_0+\cdots+\xi_{t}$ for $t\in\ZZ_+$. Assume that $\Prob(a_i\le\xi_i\le b_i)=1$ for some $\{a_i\}_{i\in\ZZ_+},\{b_i\}_{i\in\ZZ_+}\subset\R$. 
The classical Chebyshev inequality then implies that for every $\varepsilon>0$, $$ \Prob\bigr(\lvert S_{t-1}-\mathbb{E}[S_{t-1}]\rvert>\varepsilon t\bigl)\le \frac{\sum_{i=0}^{t-1}(b_i-a_i)^2}{
	\varepsilon^2t^2}.$$
In his seminal work W. Hoeffding \cite{Hoeffding-1963} has improved this result and  showed that 
$$\Prob\bigr(\lvert S_{t-1}-\mathbb{E}[S_{t-1}]\rvert>\varepsilon t\bigl)\le 2\exp\left\{-2\varepsilon^2t^2/\sum_{i=0}^{t-1}(b_i-a_i)^2\right\}.$$
Hoeffding's inequality has been widely applied   in many problems arising in probability and statistics. 
However, the independence assumption  limits its applicability  in  many situations.  This, for instance, includes problems characterized by Markovian dependence, such as Markov chain  Monte Carlo methods, time series analysis and reinforcement learning problems, see e.g.\  \cite{Fan-Jiang-Sun-2021}, \cite{Ormoneit-Glynn-2002} and \cite{Tang-2007}. Motivated by this,  Hoeffding's inequality has been extended to bounded measurable functions of a class of Markov models. 
However, to the best of our knowledge, in all the  works on this topic (see the literature review part below) a common assumption is that the underlying Markov model is irreducible. In this article, we complement these results and discuss  Hoeffding's inequality for bounded Lipschitz functions of a class of not necessarily irreducible Markov models.


\section{Main result} \label{S2}

Let $\Sp$ be a  Polish space, i.e.\
   separable completely metrizable topological space.
Denote the corresponding metric by $\mathsf{d}$. We endow $(\Sp,\mathsf{d})$ with its Borel $\sigma$-algebra
$\mathfrak{B}(\Sp)$.
  Further, let  $\T=\R_+$ or $\ZZ_+$ be
  the time parameter set, and let
$(\Omega,\mathcal{F}, \{\mathcal{F}_t\}_{t\in\T},\{\theta_t\}_{t\in\T},
\{X_t\}_{t\in\T},$ $\{\Prob_x\}_{x\in\Sp})$, denoted by $\{X_t\}_{t\in\T}$
in the sequel, be a time-homogeneous conservative strong Markov model with
 state space
$(\Sp,\mathfrak{B}(\Sp))$, in the sense of \cite{Blumenthal-Getoor-1968}. Recall that in the case when 
$\T=\ZZ_+$, $\{X_t\}_{t\in\T}$ is usually called a Markov chain, and in the case when 
$\T=\R_+$, $\{X_t\}_{t\in\T}$ is  called a Markov process.
In the latter case we also assume  that $\{X_t\}_{t\in\T}$ is progressively measurable (with respect to $\{\mathcal{F}_t\}_{t\in\T}$), i.e.\ the map $(s,\omega)\mapsto X_s(\omega)$ from $[0,t]\times\Omega$ to $\Sp$ is $\mathfrak{B}([0,t])\times\mathcal{F}_t/\mathfrak{B}(\Sp)$ measurable for all $t\geq0$. This will be in particular satisfied if $t\mapsto X_t(\omega)$ is right continuous for all $\omega\in\Omega$ (see \cite[Exercise I.6.13]{Blumenthal-Getoor-1968}).
Further, denote by
$\mathcal{P}^{t}(x,\D y)\df\mathbb{P}_{x}(X_t\in \D y)$  the transition function of $\{X_t\}_{t\in\T}$, and let 
$\mathscr{P}_1(\Sp)$ be the class of all probability measures on $\mathfrak{B}(\Sp)$ having finite first moment.
The
$\mathrm{L}^1$-Wasserstein distance on $\mathscr{P}_1(\Sp)$  is defined by
\begin{equation*}
\mathscr{W}(\upmu_1,\upmu_2)\df\inf_{\Pi\in\mathcal{C}(\upmu_1,\upmu_2)}
\int_{\Sp\times\Sp}\mathsf{d}(x,y)
\Pi(\D{x},\D{y}),
\end{equation*}
where $\mathcal{C}(\upmu_1,\upmu_2)$ is the family of couplings of
$\upmu_1(\D x)$ and $\upmu_2(\D y)$,
i.e.\ $\Pi\in\mathcal{C}(\upmu_1,\upmu_2)$ if, and only if, $\Pi(\D x,\D y)$
is a probability
measure on $\mathfrak{B}(\Sp)\times\mathfrak{B}(\Sp)$ having $\upmu_1(\D x)$ and $\upmu_2(\D y)$ as its marginals. By  Kantorovich-Rubinstein theorem it holds that
$$
\mathscr{W}(\upmu_1,\upmu_2) = \sup_{\{f\colon\mathrm{Lip}(f)\le1\}}\,
| \upmu_1(f)-\upmu_2(f)|,
$$
where the supremum is taken over all Lipschitz continuous functions
$f\colon\Sp\to\R$ with Lipschitz constant $\mathrm{Lip}(f)\le1$ and, for a probability measure $\upmu$ on $\mathfrak{B}(\Sp)$ and a measurable function $f:\Sp\to\R$, the symbol $\upmu(f)$ stands for $\int_\Sp f(x)\upmu(\D x)$, whenever the integral is well defined.

We now state the main result of this article.

\begin{theorem}\label{TM} Let $f:\Sp\to\R$ be bounded and Lipschitz continuous, and let $S_{t-1}\df\int_{[0,t)} f(X_t)\uptau(\D t)$ for $t\in\T$. Here, $\uptau(\D t)$ stands for the counting measure when $\T=\ZZ_+$ and the Lebesgue measure when $\T=\R_+$.    Assume that
	$\{X_t\}_{t\in\T}$ admits an invariant probability measure $\uppi(\D x)$ (i.e. a  measure satisfying $\int_{\Sp}\mathcal{P}^{t}(x,\D y)\uppi(\D x)=\uppi(\D y)$ for all $t\in\T$) such that \begin{equation}\label{eq:TM}\gamma\df\sup_{x\in \Sp}\int_{\T} \mathscr{W}\bigl(\mathcal{P}^{t}(x,\cdot),\uppi(\cdot)\bigr)\uptau(\D t)<\infty.\end{equation}  Then 
	 for any $\varepsilon>0$, \begin{equation}\label{eq:TM1}
	 \Prob_x\bigr(|S_{t-1}-\uppi(f)t|>t\varepsilon\bigl)\le 
	 \begin{cases}
	 2\exp\left\{\frac{-(\varepsilon t-2\mathrm{Lip}(f)\gamma)^2}{8(\mathrm{Lip}(f)\gamma+\lVert f\rVert_\infty)t}\right\}, & \T=\ZZ_+,\\[10pt]
	2 \exp\left\{\frac{-(\varepsilon t-2\mathrm{Lip}(f)\gamma)^2}{8(\mathrm{Lip}(f)\gamma+\lVert f\rVert_\infty)(t+1)}\right\}, & \T=\R_+.
	 \end{cases}\end{equation}
		\end{theorem}
	
\bigskip

According to \cite[Theorems 2.1 and 2.4]{Butkovsky-2014}	the relation in \cref{eq:TM} will hold if
	\begin{itemize}
		\item [(i)] the metric $\mathsf{d}$ is bounded (without loss of generality  by $1$)
		\item[(ii)] there is $\rho\in(0,1)$ such that for all $x,y\in\Sp$ and all $t$ large enough, $$\mathscr{W}\bigl(\mathcal{P}^{t}(x,\cdot),\mathcal{P}^{t}(y,\cdot)\bigr)\le(1-\rho)\mathsf{d}(x,y)$$
		\item[(iii)] there are  $\kappa\in\R$, measurable and bounded $\mathcal{V}\colon\Sp\to\R_+$ and concave, differentiable and  increasing to infinity function $\phi:\R_+\to\R_+$ satisfying $\phi(0)=0$, such that  \begin{equation}\label{eq:MR}\mathbb{E}_x\bigl[\mathcal{V}(X_t)\bigr]-\mathcal{V}(x)\le \kappa t-\int_{[0,t)}\mathbb{E}_x\bigl[\phi\circ\mathcal{V}(X_s)\bigr]\uptau(\D s)\end{equation}
			\item[(iv)]	there is $\epsilon\in(0,1)$ such that $$\int_{[1,\infty)}\bigl(\phi\circ\Phi^{-1}(t)\bigr)^{\epsilon-1} \uptau(\D t)<\infty,$$ where $\Phi(u)\df\int_1^u1/\phi(v)\D v$.
	\end{itemize}
Examples satisfying conditions (i)-(iv) are given in  \Cref{S5}.


\section{Literature review}\label{S3}

Hoeffding's inequality is a key tool  in the analysis of many problems arising in both probability and statistics. 
As already mentioned above, it was originally proved by W. Hoeffding \cite{Hoeffding-1963} in the context of  independent and bounded random variables. 
However,  many applied problems  require an extension of the result to the case where certain dependence of the components is involved, in particular Markovian dependence (see e.g.\ \cite{Fan-Jiang-Sun-2021}, \cite{Ormoneit-Glynn-2002} and \cite{Tang-2007}).
 Therefore, variants of the Hoeffding's inequality in the context of different types of  Markov models  have been studied recently.
 There are two main approaches to this problem: (i) based on spectral methods (see \cite{Chung-Lam-Liu-Mitzenmacher-2012}, \cite{Fan-Jiang-Sun-2021},  \cite{Leon-Perron-2004}, \cite{Lezaud-1998} \cite{Miasojedow-2014} and \cite{Rao-2019}) and (ii) based on Foster-Lyapunov inequality (see \cite{Adamczak-Bednorz-2015}, \cite{Boucher-2009}, \cite{Choi-Li-2019},  \cite{Douc-Moulines-Olsson-vanHandel-2011}, \cite{Glynn-Ormoneit-2002} and \cite{Liu-Liu-2021}). 
A common assumption  in all these works  is that the underlying Markov model is irreducible. Recall, a Markov model $\{X_t\}_{t\in\T}$ is said to be irreducible if there is a non-trivial measure $\upphi(\D x)$ on $\mathfrak{B}(\Sp)$ such that $\int_\T \mathcal{P}^t(x,B)\uptau(\D t)>0$ for all $x\in\Sp$, whenever $\upphi(B)>0$. In this article, we complement these results and obtain  Hoeffding's inequality in the case when the underlying Markov model is  not necessarily irreducible. Our result should be compared to the results obtained in \cite{Liu-Liu-2021} (see also \cite{Boucher-2009}), where Hoeffding's inequality for an irreducible Markov model $\{X_t\}_{t\in\T}$  satisfying the following  Foster-Lyapunov inequality
\begin{equation}\label{eq:LR}\mathbb{E}_x\bigl[\mathcal{V}(X_t)\bigr]-\mathcal{V}(x)\le  t-\kappa \int_{[0,t)}\mathbb{E}_x\bigl[\mathbb{1}_\mathcal{C}(X_s)\bigr]\uptau(\D s)\end{equation} has been obtained.  Compare this inequality to \cref{eq:MR}.
Here, 
  $\mathcal{V}\colon\Sp\to\R_+$ is measurable and bounded, $\kappa\in\R$ and $\mathcal{C}\in\mathfrak{B}(\Sp)$ is such that there are an atom $\alpha$ for $\{X_t\}_{t\in\T}$, $t_0\in\T$ and a non-trivial measure $\upnu(\D x)$ on $\mathfrak{B}(\Sp)$, such that $\alpha\subseteq \mathcal{C}$, $\upnu(\alpha)>0$ and $\mathcal{P}^{t_0}(x,B)\ge\upnu(B)$ for all $x\in \mathcal{C}$ and $B\in\mathfrak{B}(\Sp)$.
  Recall, a set $\alpha\in\mathfrak{B}(\Sp)$ is called an atom for $\{X_t\}_{t\in\T}$ if $\mathcal{P}^t(x,B)=\mathcal{P}^t(y,B)$ for all $x,y\in\alpha$, $t\in\T$ and $B\in\mathfrak{B}(\Sp)$. Let us remark here that this result can be slightly  generalized, i.e.\ the conclusions of \cite[Proposition 1 and Theorem 3]{Liu-Liu-2021}, and then also the main results (Hoeffding's inequality) given in \cite[Theorems 1 and 2]{Liu-Liu-2021}, remain valid by assuming \cref{eq:LR} with $\mathcal{C}$ being a petite set for $\{X_t\}_{t\in\T}$. Namely, under these assumptions in  \cite[Theorems 2.3 and 3.2]{Glynn-Meyn-1996}) it has been shown that the  solution to the corresponding stochastic Poisson equation (see \cref{eq:PE}) is uniformly bounded, which is the main step in the proof of \cite[Theorems 1 and 2]{Liu-Liu-2021}.
 Recall, a set $C\in\mathfrak{B}(\Sp)$ is said to be  petite for $\{X_t\}_{t\in\T}$ if there is probability measure $\upchi(\D t)$ on $\T$ and a non-trivial measure $\upmu(\D x)$ on $\mathfrak{B}(\Sp)$, such that $\int_\T\mathcal{P}^t(x,B)\upchi(\D t)\ge\upmu(B)$ for all $x\in C$ and $B\in\mathfrak{B}(\Sp)$. It is evident that the set $\mathcal{C}$ used in \cref{eq:LR} is petite for $\{X_t\}_{t\in\T}$.


\section{Proof of \Cref{TM}}\label{S4}
In this section, we prove \Cref{TM}. We follow and adapt the approach from \cite{Glynn-Ormoneit-2002}. By Kantorovich-Rubenstein theorem we have that $$\left|\mathbb{E}_x\bigl[f(X_t)-\uppi(f)\bigr]\right|\le \mathrm{Lip}(f) \mathscr{W}\bigl(\mathcal{P}^{t}(x,\cdot),\uppi(\cdot)\bigr).$$ Hence,  according to \cref{eq:TM} it follows  that $$\hat{f}(x)\df\int_{\T} \mathbb{E}_x\bigl[f(X_t)-\uppi(f)\bigr]\uptau(\D t)$$ is well defined and bounded. Furthermore, it clearly solves the stochastic  Poisson equation 
\begin{equation}\label{eq:PE}\mathbb{E}_x\bigl[\hat{f}(X_t)\bigr]-\hat{f}(x)=-\int_{[0,t)}\mathbb{E}_x\bigl[f(X_s)-\uppi(f)\bigr]\uptau(\D s),\end{equation} which in turn implies that $$M_t\df \hat{f}(X_t)-\hat{f}(X_0)+\int_{[0,t)}\bigl(f(X_s)-\uppi(f)\bigr)\uptau(\D s),\qquad t\in\T,$$ is a bounded martingale. 
Namely, for $s<t$ it follows that $$
\mathbb{E}_x\bigl[M_t|\mathcal{F}_s\bigr]=M_s+\mathbb{E}_{X_s}\bigl[\hat f(X_{t-s})\bigr]-\hat f(X_s)+\int_{[0,t-s)}\left(\mathbb{E}_{X_s}\bigl[f(X_u)\bigr]-\uppi(f)\right)\uptau(\D u)=M_s.
$$
By employing Markov inequality, for any $\varepsilon>0$ and $\theta\ge0$ it follows that
\begin{align*}
\Prob_x\bigr( S_{t-1}-\uppi(f)t>t\varepsilon\bigl)&\le \E^{-\theta\varepsilon t}\mathbb{E}_x\bigl[\E^{\theta(S_{t-1}-\uppi(f)t)}\bigr]\\&=\E^{-\theta\varepsilon t}\mathbb{E}_x\bigl[\E^{\theta(M_t-\hat{f}(X_t)+\hat{f}(X_0))}\bigr]\\&\le\E^{-\theta\varepsilon t+2\theta\lVert \hat f\rVert_\infty}\mathbb{E}_x\left[\exp\left\{\theta\left(\sum_{s=1}^{\lfloor t\rfloor}(M_s-M_{s-1})+M_t-M_{\lfloor t\rfloor}\right)\right\}\right].
\end{align*}
Observe that when $\T=\ZZ_+$, then $t=\lfloor t\rfloor$. Further, it clearly holds that $$|M_s-M_{s-1}|\le 2\lVert \hat f\rVert_\infty+2\lVert f\rVert_\infty\qquad \textrm{and}\qquad |M_t-M_{\lfloor t\rfloor}|\le 2\lVert \hat f\rVert_\infty+2\lVert f\rVert_\infty,$$ and
from the proof of \cite[Lemma 8.1]{Devroye-Gyorfi-Lugosi-Book-1996} it then follows that 
$$\mathbb{E}_x\left[\E^{\theta(M_s-M_{s-1})}\rvert\mathcal{F}_{s-1}\right]\le \E^{2\theta^2(\lVert \hat f\rVert_\infty+\lVert f\rVert_\infty)}\qquad \text{and} \qquad \mathbb{E}_x\left[\E^{\theta(M_t-M_{\lfloor t\rfloor})}\rvert\mathcal{F}_{\lfloor t\rfloor}\right]\le \E^{2\theta^2(\lVert \hat f\rVert_\infty+\lVert f\rVert_\infty)}.$$
Thus,
$$\mathbb{E}_x\left[\exp\left\{\theta\left(\sum_{s=1}^{\lfloor t\rfloor}(M_s-M_{s-1})+M_t-M_{\lfloor t\rfloor}\right)\right\}\right]\le \begin{cases}
\E^{2\theta^2(\lVert \hat f\rVert_\infty+\lVert f\rVert_\infty)t}, & \T=\ZZ_+,\\
\E^{2\theta^2(\lVert \hat f\rVert_\infty+\lVert f\rVert_\infty)(t+1)}, & \T=\R_+.
\end{cases}$$  
We then have 
$$
\Prob_x\bigr(S_{t-1}-\uppi(f)t>t\varepsilon\bigl)\le 
\begin{cases}

\E^{-\theta\varepsilon t+2\theta\lVert \hat f\rVert_\infty+2\theta^2(\lVert \hat f\rVert_\infty+\lVert f\rVert_\infty)t}, & \T=\ZZ_+,\\
\E^{-\theta\varepsilon t+2\theta\lVert \hat f\rVert_\infty+2\theta^2(\lVert \hat f\rVert_\infty+\lVert f\rVert_\infty)(t+1)}, & \T=\R_+.
\end{cases}$$ Analogously we conclude that  $$
\Prob_x\bigr(S_{t-1}-\uppi(f)t<-t\varepsilon\bigl)\le 
\begin{cases}

\E^{-\theta\varepsilon t+2\theta\lVert \hat f\rVert_\infty+2\theta^2(\lVert \hat f\rVert_\infty+\lVert f\rVert_\infty)t}, & \T=\ZZ_+,\\
\E^{-\theta\varepsilon t+2\theta\lVert \hat f\rVert_\infty+2\theta^2(\lVert \hat f\rVert_\infty+\lVert f\rVert_\infty)(t+1)}, & \T=\R_+.
\end{cases}$$
Finally, using $\lVert\hat f\rVert_{\infty}\le\mathrm{Lip}(f)\gamma$ (recall that $\gamma=\sup_{x\in \Sp}\int_{\T}\mathscr{W}\bigl(\mathcal{P}^{t}(x,\cdot),\uppi(\cdot)\bigr)\uptau(\D t)$) and optimizing over $\theta$ we obtain \cref{eq:TM1}.


	\section{Examples}\label{S5}
	
 In this section, we discuss several examples of non-irreducible Markov models satisfying conditions of \Cref{TM}.

 \begin{example}[Deterministic SDE] \label{EX1}{\rm  Consider the following one-dimensional (deterministic) SDE:
 		\begin{align*}\D X_t&=-|X_t|^\alpha\D t\\
 		X_0&=x\in[-1,1]
 		\end{align*} with $\alpha\in[1,2)$.
 	The SDE is well posed and it admits a unique strong solution which is a conservative strong Markov process with continuous sample paths on $\Sp=[-1,1]$ (endowed  with  the standard Euclidean metric $\mathsf{d}(x,y)=|x-y|$ and Borel $\sigma$-algebra $\mathfrak{B}([-1,1])$). When $\alpha=1$ the solution is given by $X_t=x\E^{-t}$, and when $\alpha\in(1,2)$ it is given by
 		$$X_t=
 		\frac{x}{((\alpha-1)|x|^{\alpha-1}t+1)^{1/(\alpha-1)}}.$$ 
 	Furthermore, it clearly holds that $\mathcal{P}^t(x,\D y)=\updelta_{X_t}(\D y)$, and the unique invariant probability measure of $\{X_t\}_{t\in\R_+}$ is $\updelta_{0}(\D y)$. Here,  $\updelta_{x}(\D y)$ stands for the Dirac delta measure at $x\in\Sp$.	We now have that 
 	\begin{align*} \mathscr{W}\bigl(\mathcal{P}^{t}(x,\cdot),\updelta_0(\cdot)\bigr)=|X_t|\le\begin{cases}
 	\E^{-t}, & \alpha=1,\\
 		\frac{1}{((\alpha-1)t+1)^{1/(\alpha-1)}}, & \alpha\in(1,2).
 	\end{cases} \end{align*}
 	Thus,
 	$$\sup_{x\in \Sp}\int_{0}^\infty \mathscr{W}\bigl(\mathcal{P}^{t}(x,\cdot),\updelta_0(\cdot)\bigr)\D t<\infty,$$ which shows that the condition in \cref{eq:TM} is satisfied and we can apply \Cref{TM} with any Lipschitz function $f:[-1,1]\to\R$. Observe also that $\{X_t\}_{t\in\R_+}$ is not irreducible and $\mathcal{P}^{t}(x,\D y)$ cannot converge  to $\updelta_0(\D y)$, as $t\to\infty$, in the total variation distance. \qed
 	}
 \end{example}
 
 We now give two examples of  discrete-time Markov models satisfying conditions (i)-(iv) from \Cref{S2}. We first consider an autoregressive model of order one (see e.g.\ \cite{Meyn-Tweedie-Book-2009}).
 
 \begin{example}[Autoregressive model]\label{EX2}{\rm Let $X_0$ and $\{\xi_i\}_{i\geq1}$ be random variables defined on a probability space $(\Omega,\mathcal{F},\mathbb{P})$, such that $X_0$ is independent of $\{\xi_i\}_{i\geq1}$, $\{\xi_i\}_{i\geq1}$ is an i.i.d.\  sequence,  $\mathbb{P}(X_0\in[0,1])=1$ and $\mathbb{P}(\xi_i=0)=\mathbb{P}(\xi_i=1/2)=1/2$. Define $$X_{t+1}\df\frac{1}{2}X_t+\xi_{t+1}.$$ 
 		Clearly, $\{X_t\}_{t\in\ZZ_+}$ is a Markov chain on $\Sp=[0,1]$ (endowed with the standard Euclidean metric $\mathsf{d}(x,y)=|x-y|$ and Borel $\sigma$-algebra $\mathfrak{B}([0,1])$)  with transition function $\mathcal{P}(x,\D y)=\mathbb{P}(\xi_i+x/2\in \D y)$.
 		Observe  that $\{X_t\}_{t\in\ZZ_+}$ is not irreducible. Namely, for $x\in[0,1]\cap\mathbb{Q}$ it holds that $\mathcal{P}^t(x,[0,1]\cap\mathbb{Q}^c)=0$ for all $t\ge1$, and analogously for $x\in[0,1]\cap\mathbb{Q}^c$ it holds that $\mathcal{P}^t(x,[0,1]\cap\mathbb{Q})=0$ for all $t\ge1$.
 		Next, a straightforward computation shows that $$\mathscr{W}\bigl(\mathcal{P}(x,\cdot),\mathcal{P}(y,\cdot)\bigr)\le\frac{1}{2}\mathsf{d}(x,y),$$ and since
 		\begin{equation}\label{CONT} \begin{aligned}\mathscr{W}\bigl(\mathcal{P}^{2}(x,\cdot),\mathcal{P}^{2}(y,\cdot)\bigr)&\le \inf_{\Pi\in\mathcal{C}(\mathcal{P}(x,\cdot),\mathcal{P}(y,\cdot))}
 		\int_{\Sp\times\Sp}\mathscr{W}\bigl(\mathcal{P}(u,\cdot),\mathcal{P}(v,\cdot)\bigr)
 		\Pi(\D{u},\D{v})\\&\le\frac{1}{2}\inf_{\Pi\in\mathcal{C}(\mathcal{P}(x,\cdot),\mathcal{P}(y,\cdot))}
 		\int_{\Sp\times\Sp}\mathsf{d}(u,v)
 		\Pi(\D{u},\D{v})\\
 		&=\frac{1}{2} \mathscr{W}\bigl(\mathcal{P}(x,\cdot),\mathcal{P}(y,\cdot)\bigr),\end{aligned}
 		\end{equation}
 	we conclude that $$\mathscr{W}\bigl(\mathcal{P}^t(x,\cdot),\mathcal{P}^t(y,\cdot)\bigr)\le\frac{1}{2^t}\mathsf{d}(x,y).$$	
 	Thus,	
 		 condition (ii)  from \Cref{S2} holds with  $\rho\le1/2$. Conditions (iii) and (iv) trivially hold by taking $\kappa=1$, $\mathcal{V}(x)\equiv1$ and $\phi(t)=t$. Hence, we can apply \Cref{TM} to $\{X_t\}_{t\in\ZZ_+}$ and any Lipschitz function $f:[0,1]\to\R$. Observe also that $\mathrm{Leb}(\D y)$ (on $\mathfrak{B}([0,1])$) is the (unique) invariant probability measure for $\{X_t\}_{t\in\ZZ_+}$, which is singular with respect to $\mathcal{P}^{t}(x,\D y)$ for any $t\in\ZZ_+$ and $x\in[0,1]$. Hence, $\mathcal{P}^{t}(x,\D y)$ cannot converge  to $\mathrm{Leb}(\D y)$, as $t\to\infty$, in the total variation distance. 
 		 Let us remark here that from \cite[Theorem 2.1]{Butkovsky-2014} follows that for any $\epsilon\in(0,1)$  there are $c_1(\epsilon),c_2(\epsilon)>0$, such that $$\mathscr{W}\bigl(\mathcal{P}^t(x,\cdot),\mathrm{Leb}(\cdot)\bigr)\le c_1(\epsilon) \E^{-c_2(\epsilon) t}.$$
 		 \qed }
 	 \end{example}

We now discuss a simple symmetric random walk on torus.

 \begin{example}[Random walk on torus]\label{EX3}{\rm Let $Y_0$ and $\{\xi_i\}_{i\geq1}$ be random variables defined on a probability space $(\Omega,\mathcal{F},\mathbb{P})$, such that $Y_0$ is independent of $\{\xi_i\}_{i\geq1}$, $\{\xi_i\}_{i\geq1}$ is an i.i.d.\  sequence and $\mathbb{P}(\xi_i=-1)=\mathbb{P}(\xi_i=1)=1/2$. Define $$Y_{t+1}\df Y_t+\xi_{t+1}.$$ Clearly, $\{Y_t\}_{t\in\ZZ_+}$ is a Markov model on $(\R,\mathfrak{B}(\R))$. Denote the corresponding transition  function bt $\mathcal{P}^t_Y(x,\D y).$ Next, for  $x\in\R$ let
 		$$[x]\df\bigl\{y\in\R\colon x-y\in2\pi\ZZ\bigr\},\qquad\textrm{and}\qquad
 		\mathbb{S}^1\df\bigl\{[x]\colon x\in\R\bigr\}.$$
 		Clearly,
 		$\mathbb{S}^1$ is obtained
 		by identifying the opposite
 		faces of $[0,2\pi]$. 
 		The corresponding Borel $\sigma$-algebra is denoted by $\mathfrak{B}(\mathbb{S}^1)$, which can be identified with the sub-$\sigma$-algebra of $\mathfrak{B}(\R)$  of sets of the form $\bigcup_{k\in2\pi\ZZ}\{x+k\colon x\in B\}$ for $B\in\mathfrak{B}([0,2\pi])$.
 		The covering map $\R\ni x\mapsto [x]\in\mathbb{S}^1$ is denoted by 
 		$\Pi(x)$.  The projection of  $\{Y_t\}_{t\in\ZZ_+}$, with respect to $\Pi(x)$, on the torus $\mathbb{S}^1$, denoted by   $\{X_t\}_{t\in\ZZ_+}$, is a Markov model on $(\mathbb{S}^1,\mathfrak{B}(\mathbb{S}^1))$ with transition kernel  given by
 		\begin{equation*}
 		\mathcal{P}^t(x,B)= \mathcal{P}_Y^t\bigl(z_x,\Pi^{-1}(B)\bigr)\end{equation*} for $x\in \mathbb{S}^1$, $B \in \mathfrak{B}(\mathbb{S}^1)$ and $z_x \in \Pi^{-1}(\{x\})$. Denote by 
 $\mathsf{d}(x,y)$ the arc-length metric on $\mathbb{S}^1$. It is evident that  $\mathfrak{B}(\mathbb{S}^1)$ is generated by this metric.
 It is also clear that  $\{X_t\}_{t\in\ZZ_+}$ is not irreducible. For example, for $x=[1]$ it holds that $\mathcal{P}^t(x,\Pi(\{k+2\ell\pi\colon k,l\in\ZZ\}^c))=0$ for all $t\ge1$, and analogously for $x=[\sqrt{2}]$ it holds that $\mathcal{P}^t(x,\Pi(\{k+2\ell\pi\colon k,l\in\ZZ\}))=0$ for all $t\ge1$.
		A straightforward computation shows that $$\mathscr{W}\bigl(\mathcal{P}(x,\cdot),\mathcal{P}(y,\cdot)\bigr)\le\frac{1}{2}\mathsf{d}(x,y),$$ and similarly as in \cref{CONT} 
		we conclude that $$\mathscr{W}\bigl(\mathcal{P}^t(x,\cdot),\mathcal{P}^t(y,\cdot)\bigr)\le\frac{1}{2^t}\mathsf{d}(x,y),$$	
	which is exactly
		condition (ii)  from \Cref{S2} (with  $\rho\le1/2$). As in \Cref{EX2}, conditions (iii) and (iv) trivially hold by taking $\kappa=1$, $\mathcal{V}(x)\equiv1$ and $\phi(t)=t$. Hence, we can apply \Cref{TM} to $\{X_t\}_{t\in\ZZ_+}$ and any Lipschitz function $f:\mathbb{S}^1\to\R$.
		Similarly as in  the previous example, $\mathrm{Leb}(\D y)$ (on $\mathfrak{B}(\mathbb{S}^1)$) is the (unique) invariant probability measure for $\{X_t\}_{t\in\ZZ_+}$, which is singular with respect to $\mathcal{P}^{t}(x,\D y)$ for any $t\in\ZZ_+$ and $x\in\mathbb{S}^1$. Hence, $\mathcal{P}^{t}(x,\D y)$ cannot converge  to $\mathrm{Leb}(\D y)$, as $t\to\infty$, in the total variation distance. 
		From \cite[Theorem 2.1]{Butkovsky-2014} it follows that for any $\epsilon\in(0,1)$  there are $c_1(\epsilon),c_2(\epsilon)>0$, such that $$\mathscr{W}\bigl(\mathcal{P}^t(x,\cdot),\mathrm{Leb}(\cdot)\bigr)\le c_1(\epsilon) \E^{-c_2(\epsilon) t}.$$\qed
		 }
\end{example}

\bigskip

At the end, we remark that one of  typical ways of obtaining Markov models from a given Markov model is through  a random time-change method.
Recall, a subordinator $\{S_t\}_{t\in\T_S}$ is a non-decreasing right-continuous (in the case when $\T_S=\R_+$) stochastic process
on $\R_+$ with stationary and independent increments. If $\T_S=\ZZ_+$, $\{S_t\}_{t\in\T_S}$ is a random walk; and if $\T_S=\R_+$, it is a L\'evy process.
Let now $\{X_t\}_{t\in\T}$ be  a Markov model  with
transition kernel $\mathcal{P}^t(x,\D y)$, and let
$\{S_t\}_{t\in\T_S}$ be a subordinator
independent of $\{X_t\}_{t\in\T}$. If $\T=\ZZ_+$, we assume that $\{S_t\}_{t\in\T_S}$ takes values in $\ZZ_+$.
The process $X^{S}_t\df X_{S_t}$ obtained from $\{X_t\}_{t\in\T}$ by 
a random time change through $\{S_t\}_{t\in\T_S}$, is referred to as the subordinate
process $\{X_t\}_{t\in\T}$ with subordinator $\{S_t\}_{t\in\T_S}$. 
It is easy to see that $\{X^S_t\}_{t\in\T_S}$ is again a Markov model with
transition kernel
\begin{equation*}\mathcal{P}_S^t(x,\D y)=\int_{\T_S} \mathcal{P}^s(x,\D y)\,\upmu_t(\D s),\end{equation*}
where $\upmu_t(\D s)=\mathbb{P}(S_t\in\D s)$.
It is also elementary to check that if $\uppi(\D x)$ is an invariant probability measure for
$\{X_t\}_{t\in\T}$, then it is also invariant for the subordinate process $\{X^S_t\}_{t\in\T_S}$. 
Furthermore, in \cite[Proposition 1.1]{Arapostathis-Pang-Sandric-2020} it has been shown that if 
$\mathscr{W}(\mathcal{P}^t(x,\cdot),\uppi(\cdot))\le c(x)r(t)$ for some Borel measurable $c\colon\Sp\to\R_+$ and  $r\colon\T\to\R_+$, 
then
\begin{equation*}
\mathscr{W}\bigl(\mathcal{P}_S^t(x,\cdot),\uppi(\cdot)\bigr)\le c(x) \mathbb{E}\bigl[r(S_t)\bigr].
\end{equation*}

Let us now apply this method  to Markov models from \Cref{EX1,EX2,EX3}.  Assume first that $\T_S=\ZZ_+$. In particular, this means that $\{S_t\}_{t\in\ZZ_+}$ is given as $S_t=S_{t-1}+\xi_t$, where $S_0=0$ and $\{\xi_i\}_{i\ge1}$ is a sequence of  i.i.d.\ non-negative integer-valued random variables. Assume additionally that $\Prob(\xi_i=0)=0$. This procedure is sometimes referred to as discrete subordination and it was introduced in \cite{Bendikov-Saloff-Coste-2012}. Then, in order to apply \Cref{TM} to  $\{X^S_t\}_{t\in\ZZ_+}$, it suffices to show that $\sum_{t\in\ZZ_+} \mathbb{E}[r(S_t)]<\infty.$ 
Observe that in the case of \Cref{EX1} we have that $c(x)=1$ and $$r(t)=\begin{cases}
\E^{-t}, & \alpha=1,\\
\frac{1}{((\alpha-1)t+1)^{1/(\alpha-1)}}, & \alpha\in(1,2),
\end{cases}$$ while in \Cref{EX2,EX3}, for fixed $\epsilon\in(0,1)$, $c(x)=c_1(\epsilon)$ and $r(t)=\E^{-c_2(\epsilon)t}.$ 
We now have 
\begin{align*}\sum_{t\in\ZZ_+} \mathbb{E}\bigl[r(S_t)\bigr]&=
\begin{cases}
\sum_{t\in\ZZ_+}\left(\mathbb{E}\bigl[\E^{-\xi_1}\bigr]\right)^t, & \text{\Cref{EX1} with }\alpha=1,\\
\sum_{t\in\ZZ_+}\mathbb{E}\bigl[\frac{1}{((\alpha-1)S_t+1)^{1/(\alpha-1)}}\bigr], & \text{\Cref{EX1} with }\alpha\in(1,2),\\
\sum_{t\in\ZZ_+}\left(\mathbb{E}\bigl[\E^{-c_2(\epsilon)\xi_1}\bigr]\right)^t,& \text{\Cref{EX2,EX3}},
\end{cases}\\
&\le\begin{cases}
\sum_{t\in\ZZ_+}\E^{-t}, & \text{\Cref{EX1} with }\alpha=1,\\
\sum_{t\in\ZZ_+}\frac{1}{((\alpha-1)t+1)^{1/(\alpha-1)}}, & \text{\Cref{EX1} with }\alpha\in(1,2),\\
\sum_{t\in\ZZ_+}\E^{-c_2(\epsilon)t} ,& \text{\Cref{EX2,EX3}}.
\end{cases}\end{align*}
 Thus, we can apply \Cref{TM} to $\{X^S_t\}_{t\in\ZZ_+}$ and any Lipschitz function $f:[0,1]\to\R$, $f:[-1,1]\to\R$ and, respectively, $f:\mathbb{S}^1\to\R.$ Observe also that in all three cases $\{X^S_t\}_{t\in\ZZ_+}$ is not irreducible and the corresponding transition function cannot converge to the  invariant probability measure in the total variation distance.

Let now $\T_S=\R_+$. In this case, the
 Laplace transform of $\{S_t\}_{t\in\R_+}$ takes the form
$\mathbb{E}[\E^{-uS_t}] = \E^{-t\psi(u)}$.
The characteristic (Laplace) exponent $\psi\colon(0,\infty)\to(0,\infty)$
is a Bernstein function, i.e.\ it is of class $C^\infty$ and
$(-1)^n\psi^{(n)}(u)\ge0$ for all $n\in\ZZ_+$.
It is well known that every Bernstein function admits a unique
(L\'{e}vy-Khintchine) representation 
\begin{equation*}\psi(u)=bu+\int_{(0,\infty)}(1-\E^{-uy})\,\upnu(\D y),\end{equation*}
where $b\geq0$ is the drift parameter and $\upnu(\D y)$ is a L\'{e}vy measure,
i.e.\ a Borel measure on $\mathfrak{B}\bigl((0,\infty)\bigr)$ satisfying
$\int_{(0,\infty)}(1\wedge y)\,\upnu(\D y)<\infty$.
For additional reading on  Bernstein functions we refer the reader to the
monograph \cite{Schilling-Song-Vondracek-Book-2012}.
Let now $\{S_t\}_{t\in\R_+}$ be the Poisson process (with parameter $\lambda>0$) as the simplest (non-trivial) continuous-time subordinator. Observe that in this case $b=0$ and $\upnu(\D y)=\lambda \updelta_1(\D y).$ We then have 
\begin{align*}\int_0^\infty \mathbb{E}\bigl[r(S_t)\bigr]\D t&=
\begin{cases}
\int_0^\infty\sum_{n\in\ZZ_+}\E^{-n}\frac{(\lambda t)^n}{n!}\E^{-\lambda t}\D t, & \text{\Cref{EX1} with }\alpha=1,\\
\int_0^\infty\sum_{n\in\ZZ_+}\frac{1}{((\alpha-1)n+1)^{1/(\alpha-1)}}\frac{(\lambda t)^n}{n!}\E^{-\lambda t}\D t, & \text{\Cref{EX1} with }\alpha\in(1,2),\\
\int_0^\infty\sum_{n\in\ZZ_+}\frac{1}{((\alpha-1)n+1)^{1/(\alpha-1)}}\frac{(\lambda t)^n}{n!}\E^{-\lambda t}\D t,& \text{\Cref{EX2,EX3}},
\end{cases}\\
&=\begin{cases}
\frac{\E}{\lambda(\E-1)}, & \text{\Cref{EX1} with }\alpha=1,\\
\frac{1}{\lambda}\sum_{n\in\ZZ_+}\frac{1}{((\alpha-1)n+1)^{1/(\alpha-1)}}, & \text{\Cref{EX1} with }\alpha\in(1,2),\\
\frac{\E^{c_2(\epsilon)}}{\lambda(\E^{c_2(\epsilon)}-1)} ,& \text{\Cref{EX2,EX3}}.
\end{cases}\end{align*}
 Thus, we can again apply \Cref{TM} to $\{X^S_t\}_{t\in\R_+}$ and any Lipschitz function $f:[-1,1]\to\R$, $f:[0,1]\to\R$ and, respectively, $f:\mathbb{S}^1\to\R.$ In all three cases $\{X^S_t\}_{t\in\R_+}$ is not irreducible and the corresponding transition function cannot converge to the  invariant probability measure  in the total variation distance.

So far we have considered subordinators taking values  in $\ZZ_+$ only. However, the Markov model from \Cref{EX1} can be subordinated by more general subordinators. In order to apply \Cref{TM} to such processes we again need to guarantee that \begin{equation}\label{R}\int_0^\infty\mathbb{E}\bigl[r(S_t)\bigr]\D t<\infty.\end{equation}
Note that for $\alpha\in[1,2)$ we have that $r(t)\le c(1+t)^{-\beta}$ for some $c>0$ and $\beta>1$ (possibly depending on $\alpha$). From \cite[Theorem 1.1 and Lemma 3.1]{Deng-Schilling-Song-2017} we know that if \begin{equation}\label{R2}\liminf_{s\to\infty}\frac{\psi(u)}{\log u}>0\qquad\text{and}\qquad \liminf_{s\to0}\frac{\psi(\rho u)}{\psi(u)}>1\end{equation} for some $\rho>1$, then $$\mathbb{E}\bigl[r(S_t)\bigr]\le c\left(1\wedge\psi^{-1}(1/u)\right)^{1/(\alpha-1)}.$$
Hence, \cref{R} holds if \cref{R2} and  $$\int_0^\infty\left(1\wedge\psi^{-1}(1/u)\right)^{1/(\alpha-1)}\D t<\infty$$ hold true. Typical examples of such characteristic exponents (subordinators) are given by $\psi(u)=u^\gamma$ for $\gamma\in(0,1)$ ($\gamma$-stable subordinator) and $\psi(u)=\log (1+u)$ (geometric $1$-stable subordinator).

\section*{Acknowledgements}
 Financial support through  \textit{Alexander von Humboldt Foundation} (No. HRV 1151902 HFST-E)  and \textit{Croatian Science Foundation} under project 8958 (for N.\ Sandri\'c), and the \textit{Croatian Science Foundation} under project 4197 (for S.\ \v Sebek) is gratefully acknowledged.

\bibliographystyle{alpha}
\bibliography{References}

\end{document}